\documentclass{amsart}
\usepackage{amsmath,amssymb}

\newtheorem{theorem}{Theorem}[section]

\newtheorem{lemma}[theorem]{Lemma}

\theoremstyle{definition}

\theoremstyle{remark}

\numberwithin{equation}{section}

\begin{document}

\title[A Sharp Bilinear Estimate for the Benjamin Equation]{A Sharp Bilinear Estimate for the Bourgain-type Space with Application to the Benjamin Equation}


\thanks{The first/second-named author is supported in part by NNSF of China (No.10771130)/NESEC of Canada
(No. RGPIN/261100-2003). The project is completed during the
first-named author's visit to Memorial University of Newfoundland
under the second-named author's research funding.}

\author{Wengu Chen}
\address{Institute of Applied Physics
and Computational Mathematics, P.O.Box 8009, Beijing 100088, China}
\email{chenwg@iapcm.ac.cn}

\author{Jie Xiao}
\address{Department of Mathematics and Statistics, Memorial University of Newfoundland, St. John's, NL A1C 5S7, Canada}
\email{jxiao@math.mun.ca}


\date{}



\begin{abstract}
This note shows the existence of a sharp bilinear estimate for the
Bourgain-type space and gives its application to the optimal local
well/ill-posedness of the Cauchy problem for the Benjamin equation.
\end{abstract}
\maketitle

\section{Introduction}

In the process of understanding the new function spaces and their
applications to some nonlinear evolution equations discovered by J.
Bourgain in his 1993 paper \cite{Bo}, we obtain

\begin{theorem}\label{thm1} For $\alpha,0\not=\beta,\gamma,\xi,s,b\in {\mathbb R}$ let $p(\xi)=\beta\xi^3-\alpha\xi|\xi|+\gamma\xi$ and
$X_{s,b,p}$ be the Bourgain type space -- the completion of all
$C_0^\infty(\mathbb R^2)$ functions $f$ in the norm given by
$$
\|f\|_{X_{s,b,p}}=\left(\int_{\mathbb R}\int_{\mathbb
R}(1+|\xi|^2)^{s}\big(1+|\tau-p(\xi)|^2\big)^{b}\,
|\widehat{f}(\xi,\tau)|^2 \,d\xi\,d\tau\right)^{1/2}.
$$
Then:

\item{\rm(i)} For $s>-3/4$ there exists $b\in (1/2,1)$ such
that the bilinear inequality
\begin{equation}\label{bil}
\|\partial_x(uv)\|_{X_{s,\,(b-1)+,p}}\le c
\|u\|_{X_{s,b,p}}\|v\|_{X_{s,b,p}}
\end{equation}
holds with a constant $c>0$ depending only on $s$ and $b$, where
$a\pm$ means $a\pm\epsilon$ for a sufficiently small $\epsilon>0$.

\item{\rm(ii)} For any $s\le-3/4$ and any $b\in \mathbb R$ there is
no constant $c>0$ depending only on $s$ and $b$ such that the
bilinear inequality (\ref{bil}) holds.
\end{theorem}

Here $\widehat{f}(\xi,\tau)$ is the Fourier transform in
$(\xi,\tau)$. And, it is worth pointing out that if
$(\alpha,\beta,\gamma)=(0,1,0)$ then Theorem \ref{thm1} (i) yields
Kenig-Ponce-Vega's \cite[Theorem 1.1]{KPV96}, and Theorem \ref{thm1}
(ii) implies Kenig-Ponce-Vega's \cite[Theorem 1.3]{KPV96} and
Nakanishi-Takaoka-Tsutsumi's \cite[Theorem 1 (i)]{NaTaTu}. Theorem
\ref{thm1} may seem rather specialized, but it is useful in
connection with the well/ill-posedness of the Benjamin equation.
More precisely, this theorem, together with other things, derives

\begin{theorem}\label{thm2} For $\alpha,\beta,\gamma,\xi, s,b\in {\mathbb R}$ with
$\alpha\beta\not=0$ let $p(\xi)=\beta\xi^3-\alpha\xi|\xi|+\gamma\xi$
and $H^s({\mathbb R})$ be the square Sobolev space with order $s$ --
the completion of all $C^\infty_0(\mathbb R)$ functions $f$ under
the norm
$$
\|f\|_{H^s}=\left(\int_\mathbb R
(1+|\xi|^2)^s|\widehat{f}(\xi)|^2\,d\xi\right)^\frac12.
$$
Then:

\item{\rm(i)} For $s> -3/4$ and $u_0\in H^s({\mathbb R})$ there exist
$b\in (1/2,1)$ and $T=T(\|u_0\|_{H^s})>0$ such that the Cauchy
problem for the Benjamin equation:
\begin{equation}\label{Ben1}
\left\{
\begin{array}{l}
\partial_t u - \gamma\partial_x u+\alpha {\mathcal
H}\partial^2_xu+\beta\partial^3_xu+\partial_x(u^2) =0, \quad
(x,t)\in {\mathbb R\times\mathbb
R},\\
u(x,0) =u_0(x),\quad x\in\mathbb R,
\end{array}
\right.
\end{equation}
has a unique solution $u$ in $C([0,T]; H^s({\mathbb R}))\cap
X_{s,b,p}$.

\item{\rm(ii)} For $s<-3/4$ the solution map of the above Cauchy problem is not continuous at zero, namely,
there is no $T>0$ such that the solution map of (\ref{Ben1}):
$$
u_0\in H^s({\mathbb R})\mapsto u\in C([0,T]; H^s({\mathbb R}))
$$
is continuous at zero.
\end{theorem}

Here ${\mathcal H}$ denotes the one-dimensional Hilbert transform
defined by
$$
{\mathcal H}f(x)=\lim_{\epsilon\to 0}\frac1\pi\int_{|y|>\epsilon}
{f(x-y)}y^{-1}dy,\quad x\in \mathbb R.
$$
Furthermore, four things are worth noting. The first is that the
above Benjamin equation is essentially the original Benjamin
equation \cite{B} which physically characterizes the vertical
displacement (bounded above and below by rigid horizontal planes) of
the interface between a thin layer of fluid atop and a much thicker
layer of higher density fluid (cf. \cite{Pe}) -- see
\cite{ABR}-\cite{A}-\cite{A}-\cite{CB}-\cite{LS} for the study of
existence, stability and asymptotics of solitary wave solutions of
(\ref{Ben1}). The second is that setting $\alpha\not=0$ and
$\beta=0$ in (\ref{Ben1}) generates the Benjamin-Ono equation -- see
Kenig's survey \cite{Ke} but also Ionescu-Kenig \cite{IK} and
Burq-Planchon \cite{BP} for most recent developments. The third is
that: if $(\beta,\gamma)=(-1,0)$ and $\alpha=-\nu\in (-1,0)$ then
Theorem \ref{thm2} (i) goes back to Kozono-Ogawa-Tanisaka's
\cite[Theorem 2.1]{KOT}; if $s\ge-1/8$ and $\gamma=0$ then Theorem
\ref{thm2} (i) returns to Guo-Huo's \cite[Theorem 1.1]{GH}; if
$(s,\gamma)=(0,0)$ and $\alpha\beta>0$ then Theorem \ref{thm2} (i)
is the local part of the well-posedness in Linares \cite{L}. The
fourth is that Theorem \ref{thm2} is sharp in the sense that
$s>-3/4$ and $s<-3/4$ deduce the positive and negative aspects of
the posedness of (\ref{Ben1}) respectively -- Theorem \ref{thm2}
(ii) is a new discovery and achieved via Bejenaru-Tao's argument for
\cite[Theorem 2]{BeTa} plus an example in Bourgain \cite{Bo97} and
Tzvetkov \cite{Tz} -- in the near future we will handle the
intermediate index $s=-3/4$ although it is our conjecture that
Theorem \ref{thm2} (i) can extend to this value at least in the
distributional sense just like one in Christ-Colliander-Tao's paper
\cite{ChCoTa} on the Korteweg-de Vries (KdV) equation which is
recovered from taking $(\alpha,\beta,\gamma)=(0,1,0)$ in
(\ref{Ben1}).

In order to prove Theorems \ref{thm1} and \ref{thm2} we apply T.
Tao's $[k;\,Z]$-multiplier norm method (introduced in \cite{Tao2001}
to settle some problems for the typical dispersive equations
including KdV) -- see the second and third sections of this paper,
but also the linear estimates established in
\cite{KPV93b}-\cite{Gr}-\cite{GTV} and the classical fixed point
theorem -- see the final section of this paper. For the sake of
convenience, we will use the abbreviation $<\xi>=\sqrt{1+|\xi|^2}$
for $\xi\in\mathbb R$, but also denote by $A\lesssim B$ the
statement that $A\leq CB$ holds for some large constant $C$ which
may vary from line to line and depend on various parameters;
similarly use $A\ll B$ to represent $A\leq C^{-1}B$; and use $A\sim
B$ to stand for $A\lesssim B\lesssim A$. Last but not least, we
would like to acknowledge a couple of discussions with Q.-Y. Xue
from Beijing Normal University.

\section{Fundamental Estimate for Dyadic Blocks}

From now on, for $Z$, an abelian additive group with an invariant
measure $d\xi$, and for an integer $k\geq 2$, we denote by
$\Gamma_k(Z)$ the hyperplane
$$
\Gamma_k(Z):=\{\xi=(\xi_1,\cdots,\,\xi_k)\in
Z^k:\,\xi_1+\cdots+\xi_k=0\}
$$
which is equipped with the measure
$$
\int_{\Gamma_k(Z)}f(\xi):=\int_{Z^{k-1}}f(\xi_1,\cdots,\,\xi_{k-1},\,-\xi_1
-\cdots-\xi_{k-1})d\xi_1\cdots d\xi_{k-1}.
$$
Following Tao \cite{Tao2001} we say that a function
$m:\,\Gamma_k(Z)\mapsto {\mathbb C}$ is just a $[k;\,Z]-$multiplier,
and the multiplier norm $\|m\|_{[k;\,Z]}$ is defined to be the
minimal constant $\kappa\ge 0$ such that the inequality
\begin{eqnarray*}
\Big |\int_{\Gamma_k(Z)}m(\xi)\prod_{j=1}^kf_j(\xi_j)\Big |\leq
\kappa\prod_{j=1}^k\|f_j\|_{L^2(Z)},
\end{eqnarray*}
holds for all test functions $f_j$ on $Z$. Meanwhile, we need to
review some of Tao's notations. Any summations over capitalized
variables such as $N_j,\,L_j,\,H$ are presumed to be dyadic -- that
is to say -- these variables range over numbers of the form $2^k$
for $k\in {\mathbb Z}$ (the set of all integers). If
$N_1,\,N_2,\,N_3>0$ then $N_{max}, N_{med}, N_{min}$ stand for the
maximum, median, and minimum of $N_1,\,N_2,\,N_3$ respectively, and
hence $N_{max}\ge N_{med}\ge N_{min}$. Likewise, we have
$L_{max}\geq L_{med}\geq L_{min}$ whenever $L_1,\,L_2,\,L_3>0$. More
than that, we adopt the summation conventions as follows. Any
summation of the form $L_{max}\sim\cdots$ is a sum over the three
dyadic variables $L_1,\,L_2,\,L_3\gtrsim 1$: for example,
$$
\sum_{L_{max}\sim H}:=\sum_{L_1,\,L_2,\,L_3\gtrsim 1:\,L_{max}\sim
H}.
$$
Similarly, any summation of the form $N_{max}\sim\cdots$ sum over
the three dyadic variables $N_1,\,N_2,\,N_3>0$ -- in particular --
$$
\sum_{N_{max}\sim N_{med}\sim
N}:=\sum_{N_1,\,N_2,\,N_3>0:\,N_{max}\sim N_{med}\sim N}.
$$
Finally, if $\tau,\,\xi$ and $p(\xi)$ are given with
$\tau_1+\tau_2+\tau_3=0$, then we write
$$
\lambda:=\tau-p(\xi)\quad\hbox{and}\quad\lambda_j:=\tau_j-p(\xi_j)\quad\hbox{for}\quad
j=1,\,2,\,3.
$$

In the sequel, we will establish the $[3;\,Z]-$multiplier norm
estimate for the Benjamin equation. During estimation we need the
resonance function
\begin{eqnarray}\label{resonance}
h(\xi):=p(\xi_1)+p(\xi_2)+p(\xi_3)=-\lambda_1-\lambda_2-\lambda_3,
\end{eqnarray}
which arises from what extent the spatial frequencies
$\xi_1,\,\xi_2,\,\xi_3$ share with one another. By the dyadic
decomposition of each variable $\xi_j$ or $\lambda_j$, as well as
the function $h(\xi)$, we are led to consider
\begin{eqnarray*}
\|\mathcal X\|_{[3,\,{\mathbb
R^2}]}:=\|X_{N_1,\,N_2,\,N_3;\,H;\,L_1,\,L_2,\,L_3}\|_{[3,\,{\mathbb
R}\times{\mathbb R}]},
\end{eqnarray*}
where $X_{N_1,\,N_2,\,N_3;\,H;\,L_1,\,L_2,\,L_3}$ is the multiplier
determined via
\begin{equation}\label{eqM}
X_{N_1,\,N_2,\,N_3;\,H;\,L_1,\,L_2,\,L_3}(\xi,\,\lambda):=\chi_{|h(\xi)|\sim
H}\prod_{j=1}^3\chi_{|\xi_j|\sim N_j}\chi_{|\lambda_j|\sim L_j}.
\end{equation}
From the identities
$$
\xi_1+\xi_2+\xi_3=0\quad\hbox{and}\quad\lambda_1+\lambda_2+\lambda_3+h(\xi)=0
$$
we see that $X_{N_1,\,N_2,\,N_3;\,H;\,L_1,\,L_2,\,L_3}$ vanishes
unless
\begin{equation}\label{eqU}
N_{max}\sim N_{med}\quad\hbox{and}\quad L_{max}\sim
\max\{H,\,L_{med}\}.
\end{equation}
Consequently, from (\ref{resonance}) we obtain the following
algebraic smoothing relation.

\begin{lemma}\label{le1} Let $\alpha,\beta,\gamma$ and $p(\xi)$ be the same as in Theorem \ref{thm1}. If $N_{max}\sim N_{med}\gtrsim \max\{1,\,\frac{4|\alpha|}{3|\beta|}\}$,
then
\begin{equation}\label{3sig}
\max\{|\lambda_1|, |\lambda_2|,|\lambda_3|\} \gtrsim N_1N_2N_3.
\end{equation}
\end{lemma}
\smallskip
\noindent{\bf Proof of Lemma \ref{le1}.}  Noticing
$$
p(\xi_j) =\beta \xi_j^3 -\alpha
\xi_j|\xi_j|+\gamma\xi_j\quad\hbox{for}\quad j=1,\,2,\,3,
$$
we have
\begin{eqnarray*}
h(\xi)&=& -\lambda_1- \lambda_2- \lambda_3\\
&=&p(\xi_1)+p(\xi_2)+p(\xi_3)\\
&=&3\beta\xi_1\xi_2\xi_3-\alpha(\xi_1|\xi_1|+\xi_2|\xi_2|+\xi_3|\xi_3|).
\end{eqnarray*}

Next, we simplify the last formula according to six
$(\xi_1,\xi_2)$-angle regions of $\mathbb R^2$ formed by three lines
$\xi_1=0$; $\xi_2=0$; $\xi_1+\xi_2=0$.

\smallskip

$\angle 1$: If $\xi_1\ge 0,\,\xi_2\ge 0$, then $N_3=N_{\max}$ and
hence
$$
h(\xi)=3\beta\xi_1\xi_2\xi_3-\alpha\Big(\xi_1^2+\xi_2^2-(\xi_1+\xi_2)^2\Big)
=3\beta\xi_1\xi_2\xi_3+2\alpha\xi_1\xi_2=3\beta\xi_1\xi_2
\Big(\xi_3+\frac{2\alpha}{3\beta}\Big).
$$

$\angle 2$: If $\xi_1\le 0,\,\xi_2\le 0$, then $N_3=N_{\max}$ and
hence
$$
h(\xi)=3\beta\xi_1\xi_2\xi_3-\alpha\Big(-\xi_1^2-\xi_2^2+(\xi_1+\xi_2)^2\Big)
=3\beta\xi_1\xi_2\xi_3-2\alpha\xi_1\xi_2=3\beta\xi_1\xi_2
\Big(\xi_3-\frac{2\alpha}{3\beta}\Big).
$$

$\angle 3$: If $\xi_1\ge 0,\,\xi_2\le 0,\,\xi_1+\xi_2\ge 0$, then
$N_1=N_{\max}$ and hence
$$
h(\xi)=3\beta\xi_1\xi_2\xi_3-\alpha\Big(\xi_1^2-\xi_2^2-(\xi_1+\xi_2)^2\Big)
=3\beta\xi_1\xi_2\xi_3-2\alpha\xi_2\xi_3=3\beta\xi_2\xi_3
\Big(\xi_1-\frac{2\alpha}{3\beta}\Big).
$$

$\angle 4$: If $\xi_1\ge 0,\,\xi_2\le 0,\,\xi_1+\xi_2\le 0$, then
$N_2=N_{\max}$ and hence
$$
h(\xi)=3\beta\xi_1\xi_2\xi_3-\alpha\Big(\xi_1^2-\xi_2^2+(\xi_1+\xi_2)^2\Big)
=3\beta\xi_1\xi_2\xi_3+2\alpha\xi_1\xi_3=3\beta\xi_1\xi_3
\Big(\xi_2+\frac{2\alpha}{3\beta}\Big).
$$

$\angle 5$: If $\xi_1\le 0,\,\xi_2\ge 0,\,\xi_1+\xi_2\ge 0$, then
$N_2=N_{\max}$ and hence
$$
h(\xi)=3\beta\xi_1\xi_2\xi_3-\alpha\Big(-\xi_1^2+\xi_2^2-(\xi_1+\xi_2)^2\Big)
=3\beta\xi_1\xi_2\xi_3-2\alpha\xi_1\xi_3=3\beta\xi_2\xi_3
\Big(\xi_2-\frac{2\alpha}{3\beta}\Big).
$$

$\angle 6$: If $\xi_1\le 0,\,\xi_2\ge 0,\,\xi_1+\xi_2\le 0$, then
$N_1=N_{\max}$ and hence
$$
h(\xi)=3\beta\xi_1\xi_2\xi_3-\alpha\Big(-\xi_1^2+\xi_2^2+(\xi_1+\xi_2)^2\Big)
=3\beta\xi_1\xi_2\xi_3+2\alpha\xi_2\xi_3=3\beta\xi_2\xi_3
\Big(\xi_1+\frac{2\alpha}{3\beta}\Big).
$$

As a result, we find that
$$
N_{max}\sim N_{med}\gtrsim\max\{1,\,\frac{4|\alpha|}{3|\beta|}\}
$$
implies
\begin{eqnarray*}
\max \{|\lambda_1|, |\lambda_2|, |\lambda_3|\}\ge\frac{1}{3}
(|\lambda_1 + \lambda_2+\lambda_3|) \gtrsim N_1N_2N_3,
\end{eqnarray*}
whence getting (\ref{3sig}).

Interestingly, Lemma \ref{le1} and its argument may allow us to
assume that
\begin{eqnarray}\label{eqL}
H\sim N_{\max}^2N_{\min},
\end{eqnarray}
since the multiplier in (\ref{eqM}) vanishes otherwise.

Now we are in the position to state the fundamental estimate on
dyadic blocks.

\begin{lemma}\label{block estimate} Let $\alpha,\beta,\gamma$ and $p(\xi)$ be the same as in Theorem \ref{thm1}.
Suppose that
$$
\min\{H,\,N_1,\,N_2,\,N_3,\,L_1,\,L_2,\,L_3\}>0
$$
obey (\ref{eqU}) and (\ref{eqL}). Then:

\item{\rm(i)} ((++) Coherence) $N_{max}\sim N_{min}\ \&\ L_{max}\sim H$ implies
\begin{eqnarray}\label{estimate1}
\|\mathcal X\|_{[3,\,{\mathbb R^2}]}\lesssim
L_{min}^{1/2}N_{max}^{-1/4}L_{med}^{1/4}.
\end{eqnarray}

\item{\rm(ii)} ((+-) Coherence) Anyone of the following three conditions
$$
\left\{
\begin{array}{l}
N_1\sim N_2\gg N_3\ \& \ H\sim L_3\gtrsim L_2,\,L_1;\\
N_2\sim N_3\gg N_1\ \& \ H\sim L_1\gtrsim L_2,\,L_3;\\
N_3\sim N_1\gg N_2\ \& \ H\sim L_2\gtrsim L_3,\,L_1,
\end{array}
\right.
$$
implies
\begin{eqnarray}\label{estimate2}
\|\mathcal X\|_{[3,\,{\mathbb R^2}]}\lesssim
L_{min}^{1/2}N_{max}^{-1}\Big(\min\big\{H,\,\frac{N_{max}}{N_{min}}L_{med}\big\}\Big)^{1/2}.
\end{eqnarray}

\item{\rm(iii)} In all other cases, one has
\begin{eqnarray}\label{estimate3}
\|\mathcal X\|_{[3,\,{\mathbb R^2}]}\lesssim
L_{min}^{1/2}N_{max}^{-1}\Big(\min\big\{H,\,L_{med}\big\}\Big)^{1/2}.
\end{eqnarray}
\end{lemma}

\smallskip

\noindent{\bf Proof of Lemma \ref{block estimate}.} In the high
modulation case: $L_{max}\sim L_{med}\gg H$ we have by the
elementary estimate in \cite[(37), p.861]{Tao2001},
$$
\|\mathcal X\|_{[3,\,{\mathbb R^2}]}\lesssim
L_{min}^{1/2}N_{min}^{1/2}\lesssim L_{min}^{1/2}N_{max}^{-1}
N_{min}^{1/2}N_{max}\lesssim L_{min}^{1/2}N_{max}^{-1}H^{1/2}.
$$

For the low modulation case: $L_{max}\sim H$, by symmetry we may
assume $L_1\ge L_2\ge L_3$. By \cite[Corollary 4.2]{Tao2001}, we
have
\begin{eqnarray}\label{RHD}
\|\mathcal X\|_{[3,\,{\mathbb R^2}]}&\lesssim&L_3^{1/2}\Big
|\{\xi_2:\,|\xi_2-\xi_2^0|\ll
N_{min};\,|\xi-\xi_2-\xi_3^0|\ll N_{min};\\
&&p(\xi_2)+p(\xi-\xi_2)=\tau+O(L_2)\}\Big |^{1/2}\nonumber
\end{eqnarray}
for some $\tau\in{\mathbb R},\,\xi,\,\xi_1^0,\,\xi_2^0,\,\xi_3^0$
satisfying
$$
|\xi_j^0|\sim N_j\ \hbox{for}\ j=1,\,2,\,3;\
|\xi_1^0+\xi_2^0+\xi_3^0|\ll N_{min};\ |\xi+\xi_1^0|\ll N_{min}.
$$
To estimate the right-hand side of (\ref {RHD}) we will employ the
identity
\begin{eqnarray}\label{id1}
p(\xi_2)+p(\xi-\xi_2)&=&\Big(\beta\xi_2^3-\alpha\xi_2|\xi_2|\Big)+
\Big(\beta(\xi-\xi_2)^3-\alpha(\xi-\xi_2)|\xi-\xi_2|\Big)\nonumber\\
&=&\beta\Big(\xi^3-3\xi^2\xi_2+3\xi\xi_2^2\Big)-
\alpha\Big(\xi_2|\xi_2|+(\xi-\xi_2)|\xi-\xi_2|\Big)\\
&=&q(\xi,\,\xi_2)\nonumber.
\end{eqnarray}

Now, an application of (\ref{RHD}) and (\ref{id1}) yields
\begin{eqnarray}\label{id2}
q(\xi,\,\xi_2)=\tau+O(L_2).
\end{eqnarray}
Moreover, $q(\xi,\xi_2)$ can be calculated on four angle regions of
$\mathbb R^2$ as follows.

$\angle 7$: If $\xi_2\ge 0,\,\xi-\xi_2\ge 0$, then
\begin{eqnarray*}
q(\xi,\,\xi_2)&=&\beta\Big(\xi^3-3\xi^2\xi_2+3\xi\xi_2^2\Big)-
\alpha\Big(\xi_2^2+(\xi-\xi_2)^2\Big) \\
&=&\beta\Big(\xi^3-3\xi^2\xi_2+3\xi\xi_2^2\Big)-
\alpha\Big(2\xi_2^2-2\xi\xi_2+\xi^2\Big) \\
&=&\Big(3\beta\xi-2\alpha\Big)\xi_2^2-\Big(3\beta\xi-2\alpha\Big)\xi\xi_2
+\Big(\beta\xi^3-\alpha\xi^2\Big)\\
&=&\Big(3\beta\xi-2\alpha\Big)\Big(\xi_2-\frac \xi 2\Big)^2+\frac
14\Big(\beta\xi-2\alpha\Big)\xi^2.
\end{eqnarray*}

$\angle 8$: If $\xi_2\le 0,\,\xi-\xi_2\le 0$, then
\begin{eqnarray*}
q(\xi,\,\xi_2)&=&\beta\Big(\xi^3-3\xi^2\xi_2+3\xi\xi_2^2\Big)+
\alpha\Big(\xi_2^2+(\xi-\xi_2)^2\Big) \\
&=&\Big(3\beta\xi+2\alpha\Big)\Big(\xi_2-\frac \xi 2\Big)^2+\frac
14\Big(\beta\xi+2\alpha\Big)\xi^2.
\end{eqnarray*}

$\angle 9$: If $\xi_2\ge 0,\,\xi-\xi_2\le 0$, then
\begin{eqnarray*}
q(\xi,\,\xi_2)&=&\beta\Big(\xi^3-3\xi^2\xi_2+3\xi\xi_2^2\Big)-
\alpha\Big(\xi_2^2-(\xi-\xi_2)^2\Big) \\
&=&\beta\Big(\xi^3-3\xi^2\xi_2+3\xi\xi_2^2\Big)+
\alpha\Big(\xi^2-2\xi\xi_2\Big) \\
&=&3\beta\xi\Big(\xi_2^2-(\xi+\frac{2\alpha}{3\beta})\xi_2\Big)+\Big(\beta\xi^3
+\alpha\xi^2\Big)\\
&=&3\beta\xi\Big(\xi_2-2^{-1}(\xi+\frac{2\alpha}{3\beta})\Big)^2+\Big(\frac\beta
4\xi^3-\frac{\alpha^2}{3\beta}\xi\Big).
\end{eqnarray*}

$\angle 10$: If $\xi_2\le 0,\,\xi-\xi_2\ge 0$, then
\begin{eqnarray*}
q(\xi,\,\xi_2)&=&\beta\Big(\xi^3-3\xi^2\xi_2+3\xi\xi_2^2\Big)-
\alpha\Big(-\xi_2^2+(\xi-\xi_2)^2\Big) \\
&=&3\beta\xi\Big(\xi_2-2^{-1}(\xi-\frac{2\alpha}{3\beta})\Big)^2+\Big(\frac\beta
4\xi^3-\frac{\alpha^2}{3\beta}\xi\Big).
\end{eqnarray*}

With the help of these computations, we can reach (i)-(ii)-(iii) in
Lemma \ref{block estimate}.

First, if $N_1\sim N_2\sim N_3$, then (\ref{id2}) ensures that
$\xi_2$ belongs to one interval of length
$O(L_2^{1/2}N_{\max}^{-1/2})$ no matter which one of the foregoing
four cases holds, and hence
$$
\|\mathcal X\|_{[3,\,{\mathbb R^2}]}\lesssim
L_3^{1/2}L_2^{1/4}N_{max}^{-1/4}=L_{min}^{1/2}L_{med}^{1/4}N_{max}^{-1/4},
$$
so (\ref{estimate1}) follows.

Next, if $N_2\sim N_3\gg N_1$, then
$$
|\xi_2-2^{-1}(\xi\pm\frac{2\alpha}{3\beta})|\sim N_2.
$$
From (\ref{id2}) it follows that $\xi_2$ is in one interval of
length $O(L_2N_2^{-1}N_1^{-1})$, and hence
$$
\|\mathcal X\|_{[3,\,{\mathbb R^2}]}\lesssim
L_3^{1/2}L_2^{1/2}N_1^{-1/2}N_2^{-1/2}=
L_{min}^{1/2}L_{med}^{1/2}N_{min}^{-1/2}N_{max}^{-1/2}.
$$
But $\xi_2$ is also in an interval of length $\ll N_{min}$.
Therefore (\ref{estimate2}) follows.

Last, if $N_1\sim N_2\gg N_3$ or $N_1\sim N_3\gg N_2$, then
$$
|\xi_2-2^{-1}\xi|\sim
N_1,\,|\xi_2-2^{-1}(\xi\pm\frac{2\alpha}{3\beta})|\sim N_1,
$$
and hence $\xi_2$ is in one interval of length $O(L_2N_{max}^{-2})$
no matter which one of the foregoing four cases holds. This gives
$$
\|\mathcal X\|_{[3,\,{\mathbb R^2}]}\lesssim
L_{min}^{1/2}L_{med}^{1/2}N_{max}^{-1},
$$
and consequently, (\ref{estimate3}) follows.

\section{Sharp Bilinear Estimate}

This section is devoted to verifying Theorem \ref{thm1}.

\smallskip

\noindent{\bf Proof of Theorem \ref{thm1} (i).} First of all,
Plancherel's formula tells us that proving (\ref{bil}) amounts to
showing
\begin{eqnarray}\label{norm1}
\left\|\frac{(\xi_1+\xi_2)<\xi_1>^{-s}<\xi_2>^{-s}<\xi_3>^{s}}{<\tau_1-p(\xi_1)>^{b}
<\tau_2-p(\xi_2)>^{b}<\tau_3-p(\xi_3)>^{(1-b)-}}\right\|_{[3,\,{\mathbb
R}\times{\mathbb R}]}\lesssim 1.
\end{eqnarray}
But, by the definition of $h(\xi)$ and the dyadic decomposition of
each variable $\xi_j,\,\lambda_j$ where $j=1,2,3$, we may assume
$$
|\xi_j|\sim N_j;\,|\lambda_j|\sim L_j;\,|h(\xi)|\sim H.
$$
So, using the translation invariance of the $[3;\,Z]$-multiplier
norm, we can always restrict our estimate on
$$
\min\{L_1,L_2,L_3\}\gtrsim 1\quad\hbox{and}\quad
\max\{N_1,N_2,N_3\}\gtrsim
\max\big\{1,\,\frac{4|\alpha|}{3|\beta|}\big\}.
$$
Now, the comparison principle and orthogonality from \cite[Schur's
test, p. 851]{Tao2001} reduce proving the multiplier norm estimate
(\ref{norm1}) to showing that

\begin{eqnarray}\label{id4}
\sum_{N_{max}\sim N_{med}\sim N}\sum_{L_{max}\sim L_{med}}\sum_{H\ll
L_{max}}\frac{N_3<N_3>^{s}\|\mathcal X\|_{[3;\mathbb
R^2;H]}}{<N_1>^s<N_2>^sL_1^bL_2^bL_3^{(1-b)-}}\lesssim 1
\end{eqnarray}
and
\begin{eqnarray}\label{id3}
\sum_{N_{max}\sim N_{med}\sim N}\sum_{L_1,\,L_2,\,L_3\gtrsim
1}\frac{N_3<N_3>^{s}\|\mathcal X\|_{[3;{\mathbb
R^2};{L_{\max}}]}}{<N_1>^s<N_2>^sL_1^bL_2^bL_3^{(1-b)-}}\lesssim 1
\end{eqnarray}
hold for all $N\gtrsim \max\{1,\,\frac{4|\alpha|}{3|\beta|}\}$,
where
$$
\|\mathcal{X}\|_{[3;\mathbb
R^2;H]}:=\|X_{N_1,\,N_2,\,N_3;\,H;\,L_1,\,L_2,\,L_3}\|_{[3;\,{\mathbb
R}\times{\mathbb R}]}
$$
and
$$
\|\mathcal{X}\|_{[3;\mathbb
R^2;L_{\max}]}:=\|X_{N_1,\,N_2,\,N_3;\,L_{\max};\,L_1,\,L_2,\,L_3}\|_{[3;{\mathbb
R}\times{\mathbb R}]}.
$$

Next, we verify (\ref{id3}) and (\ref{id4}). In fact, this can be
accomplished by Lemma \ref{block estimate} and the following
delicate summations in which we always fix $N\gtrsim
\max\{1,\,\frac{4|\alpha|}{3|\beta|}\}$ and consequently have
(\ref{eqL}).

We first prove (\ref{id4}). By (\ref{estimate3}) we need to verify
\begin{eqnarray}\label{id5}
&&\sum_{N_{max}\sim N_{med}\sim N}\sum_{L_{max}\sim L_{med}\gtrsim
N^2N_{min}}
\frac{N_3<N_3>^{s}L_{min}^{1/2}N_{min}^{1/2}}{<N_1>^s<N_2>^sL_1^bL_2^bL_3^{(1-b)-}}
\lesssim 1.
\end{eqnarray}
To do so, it follows from symmetry that we are required to handle
two cases:
$$
\left\{
\begin{array}{l}
 N_1\sim N_2\sim N\ \& \ N_3=N_{min};\\
 N_1\sim N_3\sim N\ \& \ N_2=N_{min}.
\end{array}
\right.
$$
Under the former case, the estimate (\ref{id5}) can be further
reduced to
\begin{eqnarray*}
&&\sum_{N_{max}\sim N_{med}\sim N}\sum_{L_{max}\sim L_{med}\gtrsim
N^2N_{min}}\frac{N^{-2s}N_{min}<N_{min}>^sL_{min}^{1/2}N_{min}^{1/2}}{L_{min}^bL_{med}^bL_{max}^{(1-b)-}}
\lesssim 1.
\end{eqnarray*}
This, after performing the $L$ summations, is reduced to
\begin{eqnarray*}
\sum_{N_{max}\sim N_{med}\sim
N}\frac{N_{min}^{3/2}<N_{min}>^s}{N^{2s}(N^{2}N_{min})^{1-}}\lesssim
1,
\end{eqnarray*}
which is true for $2+2s>0$. So, (\ref{id5}) is valid for $s>-1$.
Under the latter case, the estimate (\ref{id5}) can be reduced to
\begin{eqnarray*}
&&\sum_{N_{max}\sim N_{med}\sim N}\sum_{L_{max}\sim L_{med}\gtrsim
N^2N_{min}}\frac{NL^{1/2}_{min}N_{min}^{1/2}}{<N_{min}>^sL_{min}^bL_{med}^bL_{max}^{(1-b)-}}
\lesssim 1.
\end{eqnarray*}
However, before performing the $L$ summations, we need to pay a
little more attention to the summation of $N_{min}$. This time, we
are required to check
\begin{eqnarray*}
&&\sum_{N_{max}\sim N_{med}\sim N,\,N_{min}\leq 1}\sum_{L_{max}\sim
L_{med}\gtrsim
N^2N_{min}}\frac{NN_{min}^{1/2}}{L_{min}^{b-1/2}(L_{max})^{1-}}\\
&+&\sum_{N_{max}\sim N_{med}\sim N,\,1\leq N_{min}\leq
N}\sum_{L_{max}\sim L_{med}\gtrsim
N^2N_{min}}\frac{NN_{min}^{1/2-s}}{L_{min}^{b-1/2}(L_{max})^{1-}}\lesssim
1,
\end{eqnarray*}
which is obviously true when $s> -3/2$. Namely, (\ref{id5}) is true
for $s>-3/2$.

We then show (\ref{id3}). Under this circumstance we have
$L_{max}\sim N_{max}^2N_{min}$ and consequently consider three
matters as follows.

The first one is to handle the situation where (\ref{estimate1})
holds. In this case we have $N_1,\,N_2,\,N_3\sim N\gtrsim 1$ and so
we are required to verify
\begin{eqnarray}\label{id6}
\sum_{L_{max}\sim
N^3}\frac{N^{-s}N}{L_{min}^bL_{med}^bL_{max}^{(1-b)-}}
L_{min}^{1/2}N^{-1/4}L_{med}^{1/4}\lesssim 1.
\end{eqnarray}
Performing the $L$ summations, we are required to check
$$
\frac{N^{3/4}N^{-s}}{(N^3)^{(1-b)-}}\lesssim 1,
$$
which is true for $-3/4+s+3(1-b)>0$. So, (\ref{id6}) is true for
both $s>-3/4$ and $1/2<b<{(9+4s)}/{12}$.

The second one is to settle the case where (\ref{estimate2}) holds.
By symmetry we only need to consider two cases
$$
\left\{
\begin{array}{ll}
N\sim N_1\sim N_2\gg N_3\ \&\ H\sim L_3\gtrsim L_1,\,L_2;\\
N\sim N_2\sim N_3\gg N_1\ \&\ H\sim L_1\gtrsim L_2,\,L_3.
\end{array}
\right.
$$

For the former: $N\sim N_1\sim N_2\gg N_3\ \&\ H\sim L_3\gtrsim
L_1,\,L_2$, we are required by (\ref{estimate2}) to show
\begin{eqnarray}\label{id7}
&&\sum_{N_3\ll N}\sum_{1\lesssim L_1,\,L_2\lesssim
N^2N_3}\frac{N_3<N_3>^sL_{min}^{1/2}N^{-1}}{N^{2s}L_1^bL_2^bL_3^{(1-b)-}}
\Big(\min\Big\{N^2N_3,\,\frac{N}{N_3}L_{med}\Big\}\Big)^{1/2}\lesssim
1.
\end{eqnarray}
Splitting the left-hand side of (\ref{id7}) into two pieces $I_1$
and $I_2$ where
\begin{eqnarray*}\label{id8}
I_1=\quad\quad\sum_{N_3\leq 1}\sum_{1\lesssim L_1,\,L_2\lesssim
N^2N_3}\frac{N_3<N_3>^sL_{min}^{1/2}N^{-1}{\mathsf
M}}{N^{2s}L_1^bL_2^bL_3^{(1-b)-}}
;\\
I_2=\sum_{1<N_3\ll N}\sum_{1\lesssim L_1,\,L_2\lesssim
N^2N_3}\frac{N_3<N_3>^sL_{min}^{1/2}N^{-1}{\mathsf
M}}{N^{2s}L_1^bL_2^bL_3^{(1-b)-}}
;\\
\mathsf
M:=\Big(\min\Big\{N^2N_3,\,\frac{N}{N_3}L_{med}\Big\}\Big)^{1/2},
\end{eqnarray*}
we estimate these two terms separately. The estimate of $I_1$ goes
like this:

\begin{eqnarray}\label{id9}
I_1&\leq &\sum_{N_3\leq 1}\sum_{1\lesssim L_1,\,L_2\lesssim
N^2N_3}\frac{N_3L_{min}^{1/2}N^{-1}(\frac{N}{N_3}L_{med})^{1/2}}{N^{2s}L_{min}^{b}L_{med}^b(N^2N_3)^{(1-b)-}}
.
\end{eqnarray}
Performing the $L$ summation in (\ref{id9}), we have that
\begin{eqnarray*}
I_1&\lesssim & \sum_{N_3\leq
1}\frac{N_3^{1/2}N^{-1/2}}{N^{2s}(N^2)^{(1-b)-}N_3^{(1-b)-}}\lesssim
\frac{N^{-1/2}}{N^{2s}(N^2)^{(1-b)-}}\lesssim 1
\end{eqnarray*}
if $1/2+2s+2(1-b)>0$. Namely, $I_1\lesssim 1$ if $s>-3/4$ and
$1/2<b< 5/4+s$ are true. The estimate of $I_2$ goes like this:

\begin{eqnarray}\label{id11}
I_2&\leq &\sum_{1<N_3\ll N}\sum_{1\lesssim L_1,\,L_2\lesssim
N^2N_3}\frac{N_3N_3^sN^{-1}N^{1/2}L_{med}^{1/2}N_3^{-1/2}}{N^{2s}L_{min}^{b-1/2}L_{med}^b(N^2N_3)^{(1-b)-}}
.
\end{eqnarray}
Performing the $N_3$ summation in (\ref{id11}) we obtain that
$s-2^{-1}+b<0$ implies
\begin{eqnarray*}
I_2&\lesssim & \sum_{1\lesssim L_1,\,L_2\lesssim
N^3}\frac{1}{N^{2s}L_{min}^{b-1/2}L_{med}^{b-1/2}(N^2)^{(1-b)-}N^{1/2}}\lesssim
1.
\end{eqnarray*}
Of course, this is true if $2s+2(1-b)+1/2>0$. In other words,
$I_2\lesssim 1$ is valid for
$$
s>-3/4\quad\hbox{and}\quad\frac
12<b<\min\big\{5/4+s,\,2^{-1}-s\big\}
$$
provided $s<1/2-b<0$. If $s-1/2+b\ge 0$, then
\begin{eqnarray*}
I_2&\lesssim & \sum_{1\lesssim L_1,\,L_2\lesssim
N^3}\frac{N^{s-1/2+b}}{N^{2s}L_{min}^{b-1/2}L_{med}^{b-1/2}(N^2)^{(1-b)-}N^{1/2}}\lesssim
1
\end{eqnarray*}
is true for
$$
2s+2(1-b)+1/2>s-1/2+b.
$$
This means that $I_2\lesssim 1$ is true for $s\ge 1/2-b$ where
$1/2<b<7/8$. Combining the estimates for $I_1$ and $I_2$, we obtain
the desired estimate (\ref{id7}).

For the latter: $N\sim N_2\sim N_3\gg N_1\ \&\  H\sim L_1\gtrsim
L_2,\,L_3$, we see from (\ref{estimate2}) that (\ref{id3}) can be
established via proving
\begin{eqnarray}\label{id12}
&&\sum_{N_1\ll N}\sum_{1\lesssim L_2,\,L_3\lesssim
N^2N_1}\frac{N^{1+s}L_{min}^{1/2}N^{-1}\Big(\min\Big\{H,\,\frac
N{N_1}L_{med}\Big\}\Big)^{1/2}}{N^{s}<N_1>^sL_2^bL_3^{(1-b)-}(N^2N_1)^b}
\lesssim 1.
\end{eqnarray}
Now, writing the left-hand side of (\ref{id12}) as $J_1+J_2$ where
\begin{eqnarray*}\label{id13}
&&J_1:=\sum_{N_1\leq 1}\sum_{1\lesssim L_2,\,L_3\lesssim
N^2N_1}\frac{N^{1+s}}{N^{s}L_2^bL_3^{(1-b)-}(N^2N_1)^b}
L_{min}^{1/2}N_1^{1/2};\\
&&J_2=\sum_{1<N_1\ll N}\sum_{1\lesssim L_2,\,L_3\lesssim
N^2N_1}\frac{N^{1+s}}{N_1^sN^{s}L_2^bL_3^{(1-b)-}(N^2N_1)^b}
L_{min}^{1/2}N^{-1}N^{3/4}L_{med}^{1/4}.
\end{eqnarray*}
In $J_1$, we may assume $N_1\gtrsim N^{-2}$ -- otherwise the
summation of $L$ vanishes. Performing the summation of $L$, we get
\begin{eqnarray*}\label{id14}
J_1\lesssim \sum_{N^{-2}\lesssim N_1\leq 1}\frac{NN_1^{\frac
12-b}}{N^{2b}}\lesssim\frac{NN^{(-2)(1/2-b)}}{N^{2b}}\lesssim 1.
\end{eqnarray*}
If $1/2<b<3/4$ and $s>-3/4$ hold in $J_2$, then the summation of $L$
implies
\begin{eqnarray*}\label{id15}
J_2\lesssim \sum_{1\le N_1\ll
N}\frac{N^{3/4}N_1^{-s-b}}{N^{2b}}\lesssim
\frac{N^{3/4}N^{\max\{0,-s-b\}}}{N^{2b}}\lesssim 1.
\end{eqnarray*}
Combining the estimates for $J_1$ and $J_2$, we get the desired
estimate (\ref{id12}).

The third one is to deal with the case where (\ref{estimate3})
holds. This reduces to

\begin{eqnarray}\label{id16}
\sum_{N_{max}\sim N_{med}\sim N}\sum_{L_{max}\sim
N^2N_{min}}\frac{N_3<N_3>^{s}L_{min}^{1/2}\big(\min\{H,\,L_{med}\}\big)^{1/2}}{<N_1>^s<N_2>^sNL_1^{b}
L_2^{b}L_3^{(1-b)-}}\lesssim 1.
\end{eqnarray}
To estimate (\ref{id16}), by symmetry we need to consider two cases:
$$
\left\{
\begin{array}{l}
N_1\sim N_2\sim N\ \&\ N_3=N_{min};\\
N_1\sim N_3\sim N\ \&\ N_2=N_{min}.
\end{array}
\right.
$$

Regarding the former: $N_1\sim N_2\sim N\ \&\ N_3=N_{min}$, the
estimate (\ref{id16}) further reduces to
\begin{eqnarray*}
&&\sum_{\substack{N_1\sim N_2\sim N\\N_3\ll N}}\sum_{L_{max}\sim
N^2N_3}\frac{N_3<N_3>^s}{N^{2s}L_{min}^bL_{med}^b(N^2N_3)^{(1-b)-}}
L_{min}^{1/2}N^{-1}L_{med}^{1/2}\lesssim 1.
\end{eqnarray*}
Performing the $L$ summation, we have
\begin{eqnarray}\label{id16'}
\sum_{N_3\ll
N}\frac{N_3<N_3>^s}{N^{1+2s}(N^2)^{(1-b)-}N_3^{(1-b)-}}&=&\sum_{N_3\le
1}\frac{N_3^b}{N^{1+2s}(N^2)^{(1-b)-}}\\
&&+\sum_{1<N_3\ll
N}\frac{N_3^{s+b}}{N^{1+2s}(N^2)^{(1-b)-}}\nonumber.
\end{eqnarray}
The first term in the right-hand side of (\ref{id16'}) is bounded if
$s>-1$, while the second term in the right-hand side of
(\ref{id16'}) is less than
$\frac{N^{\max\{0,\,s+b\}}}{N^{1+2s+2(1-b)-}}$ which is bounded if
$s>-1$ and $1/2<b<3/4$. So, (\ref{id16}) is true if $s>-1$.

Regarding the latter: $N_1\sim N_3\sim N\ \& \ N_2=N_{min}$, the
estimate (\ref{id16}) can be reduced to
\begin{eqnarray}\label{id17}
&&\sum_{\substack{N_1\sim N_3\sim N\\N_2\ll N}}\sum_{L_{max}\sim
N^2N_2}\frac{N^{1+s}L_{min}^{1/2}N^{-1}\big(\min\{H,\,L_{med}\}\big)^{1/2}}{N^s<N_2>^sL_{min}^bL_{med}^bL_{max}^{(1-b)-}}
\lesssim 1.
\end{eqnarray}
Before performing the $L$ summation, we have to pay a little more
attention to the summation of $N_2$. The left-hand side of
(\ref{id17}) is now written as $J_3+J_4$ where
\begin{eqnarray*}
J_3:=\sum_{N_2\leq 1}\sum_{L_{max}\sim
N^2N_2}\frac{L_{min}^{1/2}L_{med}^{1/2}}{L_{min}^{b}L_{med}^b
L_{max}^{(1-b)-}};\\
J_4:=\sum_{1\leq N_2\leq N}\sum_{L_{max}\sim
N^2N_2}\frac{L_{min}^{1/2}L_{med}^{1/2}}{N_2^sL_{min}^{b}L_{med}^bL_{max}^{(1-b)-}}
.
\end{eqnarray*}
In $J_3$ we may assume $N_2\gtrsim N^{-2}$ -- otherwise the
summation of $L$ vanishes. Performing the summation of $L$, we get
\begin{eqnarray*}\label{id18}
J_3\lesssim \sum_{N^{-2}\lesssim N_2\leq
1}\frac{N_2^{(b-1+\epsilon)}}{N^{2(1-b-\epsilon)}}\lesssim
\frac{N^{-2(b-1+\epsilon)}}{N^{2(1-b-\epsilon)}}\lesssim 1,
\end{eqnarray*}
where $\epsilon>0$ is small enough. For $J_4$, if $s+1-b\ge 0$, then
we always have $J_4\lesssim 1$ for any $1/2<b\le 1$. If $s+1-b<0$,
we have $J_4\lesssim 1$ under $2(1-b)+s+1-b>0$. So, (\ref{id17}) is
true if $s>-3/2$ and $1/2<b<{(s+3)}/3$. This completes the proof of
the first part of Theorem \ref{thm1}.

\smallskip

\noindent{\bf Proof of Theorem \ref{thm1} (ii).} Note that
$$
\|\partial_x(uv)\|_{X_{s,b-1,p}}\lesssim\|\partial_x(uv)\|_{X_{s,(b-1)+,p}}.
$$
So it is enough to check that both $s\le-3/4$ and $b\in\mathbb R$
cannot imply
\begin{equation}\label{eqEq1}
\|\partial_x(uv)\|_{X_{s,b-1,p}}\lesssim\|u\|_{X_{s,b,p}}\|v\|_{X_{s,b,p}},
\end{equation}
which is equivalent to
\begin{equation}\label{eqEq2}
\left\|\frac{\iint_{\mathbb
R^2}\frac{f(\xi_1,\tau_1)f(\xi-\xi_1,\tau-\tau_1)(1+|\xi_1|)^{-s}(1+|\tau_1-p(\xi_1)|)^{-b}}{(1+|\xi-\xi_1|)^s(1+|\tau-\tau_1-p(\xi-\xi_1)|)^b}
\,d\xi_1d\tau_1}{\big(\frac{|\xi|(1+|\xi|)^s}{(1+|\tau-p(\xi)|)^{1-b}}\big)^{-1}}\right\|_{L^2_{\xi,\tau}(\mathbb
R^2)}\lesssim\|f\|^2_{L^2_{\xi,\tau}(\mathbb R^2)}.
\end{equation}

Case 1: $s<-3/4$. On the one hand, given a large natural number $N$
let
$$
f(\xi,\tau)=1_{A}(\xi,\tau)+1_{-A}(\xi,\tau)
$$
where $1_E$ stands for the characteristic function of a set
$E\subseteq \mathbb R^2$, and
$$
A=\big\{(\xi,\tau)\in\mathbb R^2:\ N\le\xi\le
N+N^{-1/2}\quad\hbox{and}\quad |\tau-p(\xi)|\le 1\big\};
$$
$$
-A=\{(\xi,\tau)\in\mathbb R^2:\ (-\xi,-\tau)\in A\}.
$$
See also \cite{KPV96} for the definitions of $A$ and $-A$ in the
case $p(\xi)=\xi^3$. Clearly, we have
$$
\|f\|_{L^2_{\xi,\tau}(\mathbb R^2)}\lesssim N^{-1/4}.
$$
Note that $A$ contains a rectangle with $(N,p(N))$ as a vertex, with
dimensions $10^{-2}N^{-2}\times N^{-1/2}$, and with longest side
pointing in the $(1,p'(N))$-direction where $p'(N)=3\beta
N^2-2\alpha N+\gamma$. So
$$
|(f\ast f)(\xi,\tau)|\gtrsim N^{-1/2}1_{R}(\xi,\tau)
$$
where $R$ is the rectangle centered at the origin with dimensions
$\sim N^{-2}\times N^{-1/2}$ and longest side pointing in the
$(1,p'(N))$-direction. Consequently, (\ref{eqEq2}) implies
$$
N^{-2s}N^{-1/2}N^{3(b-1)/2}N^{-1/2}N^{-1/4}\lesssim N^{-1/2}
$$
and thus $b\le1/2$. On the other hand, we also show $b>1/2$. To this
end, we apply polarization and duality to obtain that (\ref{eqEq2})
amounts to
\begin{equation*}\label{eqEq3}
\left|\iint_{\mathbb R^2}\frac{\iint_{\mathbb
R^2}\frac{f(\xi_1,\tau_1)f(\xi-\xi_1,\tau-\tau_1)<\xi_1>^{-s}<\tau_1-p(\xi_1)>^{-b}}{<\xi-\xi_1>^s<\tau-\tau_1-p(\xi-\xi_1)>^b}
\,d\xi_1d\tau_1}{\big(\frac{\xi<\xi>^sg(\xi,\tau)}{<\tau-p(\xi)>^{1-b}}\big)^{-1}}\right|\lesssim\|f\|^2_{L^2_{\xi,\tau}(\mathbb
R^2)}\|g\|_{L^2_{\xi,\tau}(\mathbb R^2)},
\end{equation*}
which is equivalent to
\begin{equation}\label{eqEq4}
\left\|\frac{\iint_{\mathbb
R^2}\frac{g(\xi,\tau)h(\xi-\xi_1,\tau-\tau_1)\xi<\xi>^{s}}{<\tau-p(\xi)>^{1-b}<\xi-\xi_1>^s<\tau-\tau_1-p(\xi-\xi_1)>^b}
\,d\xi
d\tau}{<\xi_1>^s<\tau_1-p(\xi_1)>^b}\right\|_{L^2_{\xi_1,\tau_1}(\mathbb
R^2)}\lesssim\|g\|_{L^2_{\xi,\tau}(\mathbb
R^2)}\|h\|_{L^2_{\xi,\tau}(\mathbb R^2)}.
\end{equation}
Now, if
$$
g(\xi,\tau)=1_A(\xi,\tau)\quad\hbox{and}\quad
h(\xi,\tau)=1_B(\xi,\tau),
$$
where
$$
B=\Big\{(\xi,\tau)\in\mathbb R^2:\ -N+\frac{1}{2\sqrt{N}}\le\xi\le
-N+\frac{3}{4\sqrt{N}}\quad\hbox{and}\quad |\tau-p(\xi)|\le 1\Big\}.
$$
Estimating the left-hand side of (\ref{eqEq4}) via the domain
determined by
$$
|4\tau_1-p(\xi_1)|\le4^{-1}\quad\hbox{and}\quad
2N-\frac{11}{16\sqrt{N}}\le\xi_1\le 2N-\frac{9}{16\sqrt{N}},
$$
we find
$$
N^{1/4}N^{-s}N^{-3b}\lesssim N^{-1/2}
$$
whence reaching $b>1/2$ by $s<-3/4$. Of course, we have a
contradiction.

Case 2: $s=-3/4$. From the argument for Case 1 we see that
(\ref{eqEq1}) enforces $b=1/2$. Without loss of generality we may
assume the natural numbers $N$ and $m$ are so big that $4^{m+1}\ll
N$. With this assumption and that construction in \cite[Proof of
Theorem 1 (i)]{NaTaTu} in mind we define a sequence of sets as
follows:
$$
A_j=\big\{(\xi,\tau)\in \mathbb R^2:\ N\le|\xi|\le
N+\sqrt{4^{m+1}/N}\quad\hbox{and}\quad
4^j\le|\tau-p(\xi)|<4^{j+1}\big\}
$$
for $j=0,1,...,m-1$; and $A_m$ is the union of two parallelograms
with two groups of vertices
$$
\left\{
\begin{array}{l}
\big(N,p(N)-4^m\big);\\
\big(N,p(N)-4^{m+1}\big);\\
\big(N+\sqrt{4^{m+1}/N},p(N)+p'(N)\sqrt{4^{m+1}/N}-4^m\big);\\
\big(N+\sqrt{4^{m+1}/N},p(N)+p'(N)\sqrt{4^{m+1}/N}-4^{m+1}\big),
\end{array}
\right.
$$
and
$$
\left\{
\begin{array}{l}
\big(-N,p(-N)+4^m\big);\\
\big(-N,p(-N)+4^{m+1}\big);\\
\big(-N-\sqrt{4^{m+1}/N},p(-N)-p'(-N)\sqrt{4^{m+1}/N}+4^m\big);\\
\big(-N-\sqrt{4^{m+1}/N},p(-N)-p'(-N)\sqrt{4^{m+1}/N}+4^{m+1}\big).
\end{array}
\right.
$$
Also, we set $R$ be the region comprising two parallelograms similar
to those two parallelograms making $A_m$, but with area being of a
quarter of $A_m$'s area, with two centers
$$
\Big(-\frac{7}{12}\sqrt{\frac{4^{m+1}}{N}},0\Big);\quad\Big(\frac{7}{12}\sqrt{\frac{4^{m+1}}{N}},0\Big),
$$
and with the longest sides parallel to the point
$$
\Big(\sqrt{\frac{4^{m+1}}{N}},p'(N)\sqrt{\frac{4^{m+1}}{N}}\Big).
$$
Next, given $\{a_j\}_{j=0}^m$, a finite sequence of positive numbers
let $f$ be the function on $\mathbb R^2$ decided by its Fourier
transform:
$$
\hat{f}=N\sum_{j=0}^m 4^{-j-\frac{m}{4}}a_j 1_{A_j}.
$$
Then a straightforward computation (cf. \cite[(2.6)]{NaTaTu}) gives
$$
\hat{f}\ast\hat{f}\ge N^2 a_m\Big(\sum_{j=0}^m
4^{-j-\frac{3m}{2}}a_j 1_{A_j}\Big)\ast 1_{A_m}.
$$
Applying this last inequality and noticing the following simple
facts (cf. \cite[(2.1)-(2.5)]{NaTaTu}):
$$
\hbox{Area}(A_j)\sim 4^{j+\frac{m}{2}}N^{-1/2}\quad\hbox{for}\
j=0,1,...,m-1;\quad \hbox{Area}(A_m)\sim 4^{3m/2}N^{-1/2},
$$
$$
\hbox{Area}(R)\sim 4^{3m/2}N^{-1/2};\quad
\hbox{Area}(A_j)1_R\lesssim 1_{A_j}\ast 1_{A_m}\quad\hbox{for}\
j=0,1,...,m,
$$
and
$$
j\not=k\Longrightarrow A_j\cap A_k=\emptyset,
$$
we get two groups of inequalities:
\begin{eqnarray*}
\|\partial_x(f^2)\|_{X_{-3/4,-1/2,p}}&=&\big\||\xi|(1+|\xi|)^{-3/4}(1+|\tau-p(\xi)|)^{-1/2}\hat{f}\ast\hat{f}\big\|_{L^2_{\xi,\tau}(\mathbb
R^2)}\\
&\gtrsim& \frac{\left(\iint_{R}\Big(\sum_{j=0}^m
4^{-j-3m/2}a_j(1_{A_j}\ast 1_{A_m})\Big)^2\,d\xi
d\tau\right)^{1/2}}{\big(4^{(m+1)/4}N^{3/4}a_m\big)^{-1}}\\
&\gtrsim& a_m\sum_{j=0}^m a_j
\end{eqnarray*}
and
\begin{eqnarray*}
\Big\|\frac{(1+|\tau-p(\xi)|)^{1/2}\hat{f}}{(1+|\xi|)^{3/4}}\Big\|_{L^2_{\xi,\tau}(\mathbb
R^2)}&=&\left(\iint_{\mathbb R^2}\frac{\sum_{j=0}^m
\frac{a_j^21_{A_j}}{4^{2j+m/2}}}{\Big(\frac{N^2(1+|\tau-p(\xi)|)}{(1+|\xi|)^{3/2}}\Big)^{-1}}\,d\xi d\tau\right)^{1/2}\\
&\lesssim& \Big(\sum_{j=0}^m a_j^2\Big)^{1/2},
\end{eqnarray*}
thereby reaching via (\ref{eqEq1})
\begin{equation}\label{eq5}
a_m\sum_{j=0}^m\lesssim a_j\sum_{j=0}^m a_j^2.
\end{equation}
However (\ref{eq5}) is not always true. In fact, if
$$
a_j=(1+j)^{-1}\quad\hbox{for}\quad j=0,1,...,m-1\quad\hbox{and}\quad
a_m=1,
$$
then one has a contradictory inequality:
$$
\sum_{j=1}^{m+1} j^{-1}\lesssim
\sum_{j=1}^{m+1}j^{-2}\quad\hbox{as}\quad m\to\infty.
$$
Therefore, the proof of Theorem \ref{thm2} (ii) is complete.

\section{Sharp Local Well/Ill-posedness}

This section is devoted to verifying Theorem \ref{thm2}.

\smallskip

\noindent{\bf Proof of Theorem \ref{thm2} (i).} We start with a few
notations. Denote by $W(t)$ the unitary group generating the
solution of the Cauchy problem for the linear equation
$$
\left\{
\begin{array}{l}
\partial_t v -\gamma\partial_xv+\alpha {\mathcal
H}\partial^2_xv+\beta\partial^3_xv=0, \quad (x,t)\in {\mathbb
R}\times{\mathbb R},\\
v(x,0) =v_0(x),\quad x\in\mathbb R.
\end{array}
\right.
$$
That is,
$$
v(x,t)=W(t) v_0(x) =S_t \ast v_0(x),
$$
where
$$
\widehat{S}_t=e^{itp(\xi)}\quad\hbox{or}\quad S_t(x)  =
\int_{\mathbb R} e^{i(x\xi +t p(\xi))}\,d\xi\quad\hbox{with}\quad
p(\xi)=\beta \xi^3-\alpha\xi|\xi| +\gamma\xi.
$$
Let $\psi\in C_0^\infty(\mathbb R)$ be a standard bump function such
that $\psi(t)\equiv 1$ if $|t|<1$ and $\psi(t)\equiv 0$ if $|t|>2$.
Consider the following integral equation
$$
u(x,t)= \psi(t) W(t) u_0(x) - 2\psi(\delta^{-1}t)\int_0^t W(t-t')
u(x,t')\partial_x u(x,t')\,dt'
$$
for $0<\delta<1$. Denote the right-hand side by ${\mathcal
T}(u)(x,t)$. The goal is to show that ${\mathcal T}(u)$ is a
contraction map from $Y$ to itself, where
$$
Y=\{u\in X_{s,b,p}:\,\, \|u\|_{X_{s,b,p}} \le 2c_0
 \|u_0\|_{H^s}\},
$$
where $c_0$ is the constant appeared in the following linear
estimates -- under $\alpha,0\not=\beta,\gamma\in {\mathbb R}$, one
has that for $1/2<b\le 1$,

\begin{equation}\label{ineq1}
\|\psi(\delta^{-1} t) W(t) u_0(x)\|_{X_{s,b,p}} \le c_0
\delta^{(1-2b)/2}\|u_0\|_{H^s};
\end{equation}
and for $b'+1\geq b\geq 0\geq b'>-1/2$,
\begin{equation}\label{ineq2}
\left\|\psi(\delta^{-1} t)\int_0^t W(t-t') f(x,t')\,dt'
\right\|_{X_{s,b,p}} \le c_0 \delta^{1+b'-b} \|f\|_{X_{s,b',p}}.
\end{equation}

Inequality (\ref{ineq1}) follows from Kenig-Ponce-Vega
\cite{KPV93b}, and inequality (\ref{ineq2}) follows from the
inhomogeneous linear equation version stated in \cite[Lemma
2.1]{GTV} and \cite[Lemma 1.9]{Gr}.

If $u\in Y$, then $b'$ is taken to be $b-1+\sigma$ where $\sigma>0$
is small enough to ensure that $0\geq b'>-1/2$, and hence a combined
application of (\ref{ineq1}), (\ref{ineq2}) and Theorem \ref{thm1}
(i) yields
\begin{eqnarray*}
\|{\mathcal T}(u)\|_{X_{s,b,p}}&\leq &\|\psi(t)W(t)u_0(x)\|_{X_{s,b,p}}\\
&+&2\|\psi(\delta^{-1}t)\int_0^tW(t-t')u(x,t')\partial_xu(x,t')\,dt'\|_{X_{s,b,p}}\\
&\leq &c_0\|u_0\|_{H^s}+2c_0\delta^\sigma\|u\partial_xu\|_{X_{s,b',p}}\\
&\leq &c_0\|u_0\|_{H^s}+c_0c\delta^\sigma\|u\|_{X_{s,b,p}}^2\\
&\leq &c_0\|u_0\|_{H^s}+4c_0^3c\delta^\sigma \|u_0\|_{H^s}^2.
\end{eqnarray*}
If $\delta>0$ is such a small number that
$4c_0^2c\delta^\sigma\|u_0\|_{H^s}\leq1/2$, then one has ${\mathcal
T}(u)\in Y$. Also, if $u,\,v\in Y$,
\begin{eqnarray*}
\|{\mathcal T}(u)-{\mathcal T}(v)\|_{X_{s,b,p}}&\leq
&c_0c\delta^\sigma\|u+v\|_{X_{s,b,p}}\|u-v\|_{X_{s,b,p}}\\
&\leq
&4c_0^2c\delta^\sigma\|u_0\|_{H^s}\|u-v\|_{X_{s,b,p}}\\
&\leq&2^{-1}\|u-v\|_{X_{s,b,p}}.
\end{eqnarray*}
Therefore ${\mathcal T}$ is a contraction mapping on $Y$. By the
classical Banach fixed point theorem, there exists a unique solution
$u\in Y$ such that
$$
u(x,t)= \psi(t) W(t) u_0(x) - 2\psi(\delta^{-1}t)\int_0^t W(t-t')
u(x,t')
\partial_x u(t')\,dt'.
$$
Choosing $T=2^{-1}\delta$, we have that $t\in [0,T]$ implies
$\psi(t)=1$ and so that $u(x,t)$ solves the integral equation
associated to the Cauchy problem (\ref{Ben1}).

Next, we verify the persistence property $u\in C([0,T],H^s(\mathbb
R))$ and the continuous dependence of the solution upon the data.
Clearly, the former follows directly from \cite[Corollary 2.1]{Tao}
which gives
$$
\sup_{t\in [0,T]}\|u(\cdot,t)\|_{H^s}\lesssim\|u\|_{X_{s,b,p}}.
$$
As to the latter, we apply (\ref{ineq1})-(\ref{ineq2}),
\cite[Corollary 2.1\ \&\ Lemma 2.64]{Tao} (with $\eta$ being a
Schwartz function on $\mathbb R$, e.g., $\eta=\psi$ as above) and
Theorem \ref{thm1} (i) to obtain that if $0\le t_0<t\le T$ and
$t-t_0\le\Delta t$ then
\begin{eqnarray*}
\|u(\cdot,t)-u(\cdot,t_0)\|_{H^s}&\lesssim&\|W(t-t_0)u(\cdot,t_0)-u(\cdot,t_0)\|_{H^s}\\
&&\quad+\left\|\int_{t_0}^tW(t-t')\eta^2\Big(\frac{t'-t_0}{\Delta
t}\Big)\partial_x\big(u^2(\cdot,t')\big)\,dt'\right\|_{H^s}\\
&\lesssim&\|W(t-t_0)u(\cdot,t_0)-u(\cdot,t_0)\|_{H^s}\\
&&\quad+\left\|\int_{t_0}^tW(t-t')\eta^2\Big(\frac{t'-t_0}{\Delta
t}\Big)\partial_x\big(u^2(\cdot,t')\big)\,dt'\right\|_{X_{s,b,p}}\\
&\lesssim&\|W(t-t_0)u(\cdot,t_0)-u(\cdot,t_0)\|_{H^s}\\
&&\quad+\left\|\eta^2\Big(\frac{t'-t_0}{\Delta t}\Big)\partial_x\big(u^2(\cdot,t')\big)\right\|_{X_{s,b-1,p}}\\
&\lesssim&\|W(t-t_0)u(\cdot,t_0)-u(\cdot,t_0)\|_{H^s}\\
&&\quad+(\Delta t)^{0+}\left\|\partial_x\big(u^2(\cdot,t')\big)\right\|_{X_{s,(b-1)+,p}}\\
&\lesssim&\|W(t-t_0)u(\cdot,t_0)-u(\cdot,t_0)\|_{H^s}+(\Delta
t)^{0+}\|u\|^2_{X_{s,b,p}}\\
&=&o(1)\quad\hbox{as}\quad \Delta t\to 0,
\end{eqnarray*}
giving the persistence property.

To close the argument we need to demonstrate that the uniqueness of
the solution to (\ref{Ben1}). To this end, for $\tau>0$ let
$$
\|u\|_{X_{s,b,p}^\tau}=\inf\big\{\|v\|_{X_{s,b,p}}:\ v\in
X_{s,b,p}\quad\hbox{and}\quad
v(\cdot,t)=u(\cdot,t)\quad\hbox{for}\quad t\in [0,\tau]\big\}.
$$
Clearly, if $\|u_1-u_2\|_{X_{s,b,p}^\tau}=0$, then
$u_1(\cdot,t)=u_2(\cdot,t)$ in $H^s(\mathbb R)$ for $t\in [0,\tau]$.
Suppose now that $u_1$ is the solution on $[0,T]$ obtained by the
fixed point theorem as above and $u_2$ is a solution of the integral
equation associated to (\ref{Ben1}) with the same initial data
$u_0$. Without loss of generality, we may assume that there is a
constant $M>1$ such that
$$
\max\{\|u_1\|_{X_{s,b,p}},\|\eta u_2\|_{X_{s,b,p}}\}\le M
$$
where $\eta$ is the above-appeared bump function on $\mathbb R$, and
also $T\in (0,1)$. Then for $0<T^\ast<T$ one has
$$
\eta(t) u_2(x,t)=\eta(t)W(t)u_0(x)-\eta(t)\int_0^t
W(t-t')\eta\Big(\frac{t'}{T^\ast}\Big)\eta^2(t')\partial_x(u_2^2(x,t'))\,dt'
$$
where $(x,t)\in \mathbb R\times[0,T^\ast]$.

For any $\epsilon>0$ there exists $w\in X_{s,b,p}$ such that $t\in
[0,T^\ast]$ implies
$$
w(x,t)=u_1(x,t)-\eta(t)u_2(x,t)\quad\hbox{and}\quad
\|w\|_{X_{s,b,p}}\le\|u_1(x,t)-\eta(t)
u_2(x,t)\|_{X^{T^\ast}_{s,b,p}}+\epsilon.
$$
So if
$$
\tilde{w}(x,t)=-\eta(t)\int_0^t
W(t-t')\eta\Big(\frac{t'}{T^\ast}\Big)\partial_x\big(w(x,t')u_1(x,t')+\eta(t')w(x,t')u_2(x,t')\big)\,dt',
$$
then
$$
\tilde{w}(x,t)={w}(x,t)=u_1(x,t)-\eta(t)u_2(x,t)\quad\hbox{for}\quad
t\in [0,T^\ast].
$$
According to the linear estimates (\ref{ineq1})-(\ref{ineq2}) and
Theorem \ref{thm1} (i), for $-1/2<b'=b-1+\sigma\le 0$ (with
$0<\sigma$ being small enough) we have
\begin{eqnarray*}
\|u_1(x,t)&-&\eta(t)u_2(x,t)\|_{X^{T^\ast}_{s,b,p}}\\
&\le&\|\tilde{w}\|_{X_{s,b,p}}\\
&\lesssim&(T^\ast)^\sigma\big\|\eta(t'/T^\ast)\partial_x\big(w(x,t')u_1(x,t')+\eta(t')w(x,t')u_2(x,t')\|_{X_{s,b',p}}\\
&\lesssim&(T^\ast)^\sigma\big(\|w\|_{X_{s,b,p}}\|u_1\|_{X_{s,b,p}}+\|w\|_{X_{s,b,p}}
\|\eta u_2\|_{X_{s,b,p}}\big)\\
&\lesssim&M(T^\ast)^\sigma\|w\|_{X_{s,b,p}},
\end{eqnarray*}
which produces a constant $c_1>0$ such that
$$
\|u_1(x,t)-\eta(t)u_2(x,t)\|_{X^{T^\ast}_{s,b,p}}\le
c_1M(T^\ast)^\sigma\|w\|_{X_{s,b,p}}.
$$
If $T^\ast\le(2c_1M)^{-1/\sigma}$ then
$$
\|u_1(x,t)-\eta(t)u_2(x,t)\|_{X_{s,b,p}^{T^\ast}}\le
2^{-1}\|w\|_{X_{s,b,p}}\le2^{-1}(\|u_1(x,t)-\eta(t)u_2(x,t)\|_{X_{s,b,p}^{T^\ast}}+\epsilon)
$$
and hence
$$
\|u_1(x,t)-\eta(t)u_2(x,t)\|_{X_{s,b,p}^{T^\ast}}\le 2\epsilon.
$$
Since $\epsilon>0$ is arbitrary, the last inequality yields that
$u_1(\cdot,t)=u_2(\cdot,t)$ for all $t\in [0,T^\ast]$. Continuing
this process, we achieve the uniqueness assertion on $[0,T]$.

\smallskip

\noindent{\bf Proof of Theorem \ref{thm2} (ii).} Under $s<-3/4$ for
contradiction we assume that the solution map

$$
u_0\in H^s({\mathbb
R})\mapsto u\in C([0,T]; H^s({\mathbb R}))
$$
is continuous at zero. According to Bejenaru-Tao's \cite[Theorem 3
\& Proposition 1]{BeTa} -- a general principle for well-posedness,
we must have
\begin{equation*}\label{normI}
\sup_{t\in
[0,T]}\|A_3(f)\|_{H^s}\lesssim\|f\|_{H^s}^3\quad\hbox{for\ all}\quad
f\in H^s(\mathbb R),
\end{equation*}
where
$$
A_3(f)(x,t)=\int_{\mathbb R}e^{ix\xi}\left(\int_0^t
\frac{(i\xi)\Big(\int_{\mathbb
R}\widehat{A_1(f)}(\xi',t')\widehat{A_2(f)}(\xi-\xi',t')\,d\xi'\Big)}{e^{i(t'-t)p(\xi)}}\,dt'\right)\,d\xi;
$$
$$
A_2(f)(x,t)=\int_{\mathbb R}\frac{\left(\int_{\mathbb
R}\frac{\Big(\int_0^t
e^{it'\big(p(\xi_1)+p(\xi_2)-p(\xi_1+\xi_2)\big)}\,dt'\Big)}{\Big(e^{itp(\xi_1+\xi_2)}\big(i(\xi_1+\xi_2)\big)\hat{f}(\xi_1)\hat{f}(\xi_2)\Big)^{-1}}\,d\xi_1\right)}{e^{-ix(\xi_1+\xi_2)}}\,d\xi_2;
$$
$$
A_1(f)(x,t)=\int_{\mathbb R}e^{itp(\xi)+ix\xi}\hat{f}(\xi)\,d\xi.
$$
In the definition of $A_3$, the Fourier transform is taken over the
spatial variable.

Motivated by the selection of a test function in \cite{Bo97} and
\cite{Tz} we choose an $H^s(\mathbb R)$-function $f$ with
$$
\|f\|_{H^s}\sim
1\quad\hbox{and}\quad\hat{f}(\xi)=r^{-1/2}N^{-s}1_{[-r,r]}(|\xi|-N),
$$
where $r=(\sqrt{N}\log N)^{-1}$, $N>0$ is sufficiently large, and
$1_{E}$ stands for the characteristic function of a set
$E\subseteq\mathbb R$.

The key issue is to control $\|A_3(f)\|_{H^s}$ from below. To
proceed, we make the following estimates:
$$
A_1(f)(x,t)\sim r^{-1/2}N^{-s}\int_{|\xi\pm
N|<r}e^{itp(\xi)+ix\xi}\,d\xi
$$
and
$$
A_2(f)(x,t)\sim F_1(x,t)-F_2(x,t)
$$
where
$$
F_1(x,t)=r^{-1}N^{-2s}\iint_{\max_{j=1,2}|\xi_j\pm
N|<r}\frac{(\xi_1+\xi_2)e^{ix(\xi_1+\xi_2)+it\big(p(\xi_1)+p(\xi_2)\big)}}{p(\xi_1)+p(\xi_2)-p(\xi_1+\xi_2)}\,d\xi_1d\xi_2
$$
and
$$
F_2(x,t)=r^{-1}N^{-2s}\iint_{\max_{j=1,2}|\xi_j\pm
N|<r}\frac{(\xi_1+\xi_2)e^{ix(\xi_1+\xi_2)+it(p(\xi_1+\xi_2))}}{p(\xi_1)+p(\xi_2)-p(\xi_1+\xi_2)}\,d\xi_1d\xi_2.
$$

The contribution of $F_1$ to $A_3(f)$ is comparable with
\begin{equation}\label{com}
r^{-3/2}N^{-3s}\iiint_{\max_{j=1,2,3}|\xi_j\pm
N|<r}\frac{Q_1(\xi_1,\xi_2,\xi_3)Q_2(\xi_1,\xi_2,\xi_3)}{e^{-ix(\xi_1+\xi_2+\xi_3)-it(p(\xi_1+\xi_2+\xi_3))}}\,d\xi_1d\xi_2d\xi_3,
\end{equation}
where
$$
Q_1(\xi_1,\xi_2,\xi_3):=\frac{(\xi_1+\xi_2+\xi_3)(\xi_2+\xi_3)}{p(\xi_2)+p(\xi_3)-p(\xi_2+\xi_3)}
$$
and
$$
Q_2(\xi_1,\xi_2,\xi_3):=\frac{e^{it\big(p(\xi_1)+p(\xi_2)+p(\xi_3)-p(\xi_1+\xi_2+\xi_3)\big)}-1}{
p(\xi_1)+p(\xi_2)+p(\xi_3)-p(\xi_1+\xi_2+\xi_3)}.
$$
Setting
$$
\theta=p(\xi_1)+p(\xi_2)+p(\xi_3)-p(\xi_1+\xi_2+\xi_3)\quad\hbox{and}\quad
-\xi_4=\xi_1+\xi_2+\xi_3,
$$
we employ $p(\xi)=\beta\xi^3-\alpha\xi|\xi|+\gamma\xi$ to get
$$
\theta=\beta(\xi_1^3+\xi_2^3+\xi_3^3+\xi_4^3)-\alpha(\xi_1|\xi_1|+\xi_2|\xi_2|+\xi_3|\xi_3|+\xi_4|\xi_4|).
$$
By symmetry we may assume that
$|\xi_1|\ge|\xi_2|\ge|\xi_3|\ge|\xi_4|$ and further $\xi_1\ge 0$
thanks to $\sum_{j=1}^4\xi_j=0$, and consequently consider two cases
according to the signs of $\xi_1$, $\xi_2$, $\xi_3$ and $\xi_4$.

Case 1: (+,-,-,-). In this case we have
\begin{eqnarray*}
\theta&=&3\beta(\xi_1+\xi_4)(\xi_2+\xi_4)(\xi_3+\xi_4)-\alpha(\xi_1^2-\xi_2^2-\xi_3^2-\xi_4^2)\\
&=&3\beta(\xi_1+\xi_4)(\xi_2+\xi_4)(\xi_3+\xi_4)-\alpha\big((\xi_1-\xi_2)(\xi_1+\xi_2)-(\xi_3+\xi_4)^2+2\xi_3\xi_4\big)\\
&=&3\beta(\xi_1+\xi_4)(\xi_2+\xi_4)(\xi_3+\xi_4)-\alpha\big((\xi_3+\xi_4)(\xi_2-\xi_1-\xi_3-\xi_4)+2\xi_3\xi_4\big)\\
&=&3\beta(\xi_1+\xi_4)(\xi_2+\xi_4)(\xi_3+\xi_4)-\alpha\big((\xi_3+\xi_4)(2\xi_2)+2\xi_3\xi_4\big),
\end{eqnarray*}
whence finding

$$
\theta\sim
\left\{
\begin{array}{l}
\xi_1\xi_2\xi_3\quad\hbox{if}\quad |\xi_4|\ll
|\xi_1|,\\
\xi_2^2\xi_3\quad\hbox{if}\quad |\xi_4|\sim|\xi_1|.
\end{array}
\right.
$$

Case 2: (+,-,-,+). In this case we have
\begin{eqnarray*}
\theta&=&\beta(\xi_1^3+\xi_2^3+\xi_3^3+\xi_4^3)-\alpha(\xi_1^2-\xi_2^2-\xi_3^2-\xi_4^2)\\
&=&3\beta(\xi_1+\xi_4)(\xi_2+\xi_4)(\xi_3+\xi_4)-2\alpha(\xi_3+\xi_4)(\xi_2+\xi_4)\\
&=&3\beta(\xi_2+\xi_4)(\xi_3+\xi_4)\Big(\xi_1+\xi_4-\frac{2\alpha}{3\beta}\Big)\\
&\sim&(\xi_2+\xi_4)(\xi_3+\xi_4)\xi_1\\
&\sim&(\xi_1+\xi_3)(\xi_1+\xi_2)\xi_1.
\end{eqnarray*}
Thus
$$
|\theta|\sim N^3\quad\hbox{or}\quad |\theta|\lesssim r^2N\sim (\log
N)^{-2}.
$$
This tells us that the major contribution to (\ref{com}) is obtained
via
$$
G_1(x,t)=\frac{\iiint_{\max_{j=1,2,3}\{|\xi_j\pm N|\}<r,\
|\theta|\lesssim
r^2N}\frac{Q_1(\xi_1,\xi_2,\xi_3)}{e^{-i\big(x(\xi_1+\xi_2+\xi_3)+
tp(\xi_1+\xi_2+\xi_3)\big)}}\,d\xi_1d\xi_2d\xi_3}{r^{3/2}N^{3s}}
$$
with
$$
\|G_1\|_{H^s}\sim rN^{-2s-1}\sim {N^{-2s-3/2}}{(\log N)^{-1}}.
$$

On the other hand, the contribution of $F_2$ to $A_3(f)$ is
comparable with
$$
G_2(x,t)=\frac{\iiint_{\max_{j=1,2,3}|\xi_j\pm
N|<r}\frac{Q_1(\xi_1,\xi_2,\xi_3)Q_3(\xi_1,\xi_2,\xi_3)}{e^{-i\big(x(\xi_1+\xi_2+\xi_3)+
tp(\xi_1+\xi_2+\xi_3)\big)}}\,d\xi_1d\xi_2d\xi_3}{r^{3/2}N^{3s}},
$$
where
$$
Q_3(\xi_1,\xi_2,\xi_3):=\frac{e^{it\big(p(\xi_1)+p(\xi_2+\xi_3)-p(\xi_1+\xi_2+\xi_3)\big)}-1}{p(\xi_1)+p(\xi_2+\xi_3)-p(\xi_1+\xi_2+\xi_3)}
$$
and
\begin{eqnarray*}
\|G_2\|_{H^s}&\lesssim& \frac{\left\|\iiint_{\max_{j=1,2,3}|\xi_j\pm
N|<r}\frac{e^{ix(\xi_1+\xi_2+\xi_3)}}{|\xi_2+\xi_3|+N^{-2}}\,d\xi_1d\xi_2d\xi_3\right\|_{L^2_{x}(\mathbb
R)}}{r^{3/2}N^{2s+3}}\\
&\lesssim& r^{-1}N^{-2s-3}\iint_{\max_{j=2,3}|\xi_j\pm N|<r}(|\xi_2+\xi_3|+N^{-2})^{-1}\,d\xi_2d\xi_3\\
&\lesssim& N^{-2s-3}\log N.
\end{eqnarray*}
Consequently, we get
$$
\frac{N^{-2s-3/2}}{\log N}\left(1-\Big(\frac{\log
N}{N^{3/4}}\Big)^2\right)\lesssim\|G_1\|_{H^s}-\|G_2\|_{H^s}\lesssim
\|A_3(f)\|_{H^s}\lesssim 1
$$
whence deriving $s\ge -3/4$ (via letting $N\to\infty$) -- a
contradiction to $s<-3/4$. This completes the proof of Theorem
\ref{thm2} (ii).

\bibliographystyle{amsplain}

\end{document}